\documentclass[journal]{IEEEtran}
\usepackage{cite}

\ifCLASSINFOpdf
  \usepackage[pdftex]{graphicx}
  \usepackage{epstopdf}
  \usepackage{subfigure}
  \graphicspath{{./}{./figures/}}
\else
\fi
\usepackage{xcolor}

\usepackage{amsmath}
\usepackage{amssymb}
\newcommand*{\qed}{\hfill\ensuremath{\square}}%

\newtheorem{theorem}{Theorem}[section]
\newtheorem{corollary}{Corollary}[theorem]
\newtheorem{lemma}[theorem]{Lemma}

\newtheorem{assum}[theorem]{Assumption}

\newcommand{\norme}[1]{\left\Vert #1\right\Vert}

\begin{document}
%
\title{Boundary Control of a Nonhomogeneous Flexible Wing with Bounded Input Disturbances}
%
%
%

\author{Hugo~Lhachemi,~\IEEEmembership{Student~Member,~IEEE,}
        David~Saussi{\'e},~\IEEEmembership{Member,~IEEE,}
        and~Guchuan~Zhu,~\IEEEmembership{Senior~Member,~IEEE}
\thanks{The authors are with the Department of Electrical Engineering, {\'E}cole Polytechnique de Montr{\'e}al, Montr{\'e}al, QC, H3T-1J4 Canada e-mail: \{hugo.lhachemi,d.saussie,guchuan.zhu\}@polymtl.ca.}
}

%
%

\markboth{Manuscript submitted to IEEE Transaction on Automatic Control}%
{Lhachemi \MakeLowercase{\textit{et al.}}: Boundary Control for Flutter Suppression of a Nonhomogeneous Flexible Wing}
%



\maketitle

\begin{abstract}
This note deals with the boundary control problem of a nonhomogeneous flexible wing evolving under unsteady aerodynamic loads. The wing is actuated at its tip by flaps and is modeled by a distributed parameter system consisting of two coupled partial differential equations. Based on the proposed boundary control law, the well-posedness of the underlying Cauchy problem is first investigated by resorting to the semigroup theory. Then, a Lyapunov-based approach is employed to assess the stability of the closed-loop system in the presence of bounded input disturbances.
\end{abstract}

\begin{IEEEkeywords}
Distributed parameter systems, Flexible structures, Aerospace systems, Lyapunov stability.
\end{IEEEkeywords}

%
\IEEEpeerreviewmaketitle

\section{Introduction}\label{sec: Introduction}
%
%
%
%
Partial Differential Equation (PDE)-based control of flexible structures has attracted much attention in the last decades. For instance, control of Euler-Bernoulli beams is one of the most investigated problems for which different design strategies have been applied, including, e.g., backstepping control~\cite{Krstic2008}, Lyapunov method~\cite{Queiroz2012,Luo2012}, passivity-based control~\cite{henikl2016}, flatness method~\cite{Aoustin1997,Meurer2012}, spectral analysis~\cite{Curtain2012,Luo2012}, and optimal control~\cite{Bensoussan2007,Curtain2012}. It is also reported in recent literature that PDE-based control of flexible aircraft wing modeled by coupled beam and string equations, describing bending and twisting displacements, has been applied to conventional aircraft or UAV flapping wings~\cite{Paranjape2013,Bialy2016,He2017,2017arXiv170310182L}. 

The aforementioned work has considered homogeneous wings for which the parameters, such as mass and or rigidity, are supposed to be constant along the wingspan. Obviously, this consideration is not truly representative for real-life aircraft wings that should be more accurately modeled as nonhomogeneous structures. The stabilization of the bending dynamics of nonhomogeneous beams has been investigated in~\cite{guo2002riesz,chentouf2006stabilization,chentouf2007optimal,chen2015exponential}. Nevertheless, the control of the coupled bending and twisting dynamics of a nonhomogeneous wing is more challenging, because of the inherent difficulty to establish the well-posedness of such complex systems. Moreover, for stability assessment, the method of spectral analysis as proposed in~\cite{guo2002riesz,chentouf2006stabilization,chentouf2007optimal,chen2015exponential} might not be applicable due to the difficulty to find the closed form eigenfunctions of the considered coupled PDEs. As a continuous development of the work presented 
 in~\cite{2017arXiv170310182L}, this note addresses the problem of boundary stabilization of a nonhomogeneous wing under unsteady aerodynamic loads with actuators located at the wing tip. The wing is modeled as a distributed system composed of two coupled PDEs with asymmetric structures and nonconstant coefficients, describing the bending and twisting dynamics along the wingspan~\cite{Bialy2016,Zhang2005,Ziabar2010}.

Compared to the work presented in~\cite{2017arXiv170310182L}, the contribution of this note is twofold. First, it is shown that the control law proposed in~\cite{2017arXiv170310182L} for a homogeneous flexible wing applies also to the stabilization of the considered nonhomogeneous structure, which is a more practically relevant problem. Second, the impact of input disturbances on the stability properties is investigated. Specifically, the problem is formulated under an abstract form, allowing the application of the semigroup theory~\cite{Curtain2012,Pazy2012}. In particular, it is shown that the closed-loop system with the proposed boundary control is well-posed. Then, a Lyapunov-based stability analysis is performed, which shows that under certain structural constraints of the wing physical parameters, the underlying $C_{0}$-semigroup is exponentially stable. Finally, the impact of bounded input disturbances on the closed-loop system stability is evaluated.

The remainder of the note is organized as follows. Notations and preliminaries are presented in Section~\ref{sec: notations}. The wing model and the associated abstract form are introduced in Section~\ref{sec: model}. The well-posedness of the problem is tackled in the framework of semigroup theory in Section~\ref{sec: well-posedness}. A Lyapunov-based stability analysis is carried out in Section~\ref{sec: stability}. The temporal behavior of the closed-loop system is evaluated in Section~\ref{sec: simulations}, followed by some concluding remarks in Section~\ref{sec: conclusion}.

\section{Notations and preliminaries}\label{sec: notations}
The sets of real, non-negative real, positive real, and complex numbers are denoted by $\mathbb{R}$, $\mathbb{R}_+$, $\mathbb{R}_+^*$, and $\mathbb{C}$, respectively. For any given Lebesgue measurable function $f$ from $(0,l)$ to $\mathbb{R}$, the essential supremum and the essential infimum of $f$ are defined respectively by
\begin{equation*}
\overline{f} \triangleq \inf \left\{ M\in\mathbb{R} \; : \; \lambda(f^{-1}(M,+\infty)) = 0 \right\}  ,
\end{equation*}  
\begin{equation*}
\underline{f} \triangleq \sup \left\{ m\in\mathbb{R} \; : \; \lambda(f^{-1}(-\infty,m)) = 0 \right\}  ,
\end{equation*}  
where $\lambda$ stands for the Lebesgue measure. The set of Lebesgue measurable functions $f$ from $(0,l)$ to $\mathbb{R}$ which are essentially bounded, i.e., for which $\overline{|f|} < \infty$, is denoted by $L^\infty(0,l)$ and is endowed with the norm $\norme{f}_{L^\infty(0,l)} = \overline{|f|}$. For the set of continuous functions over a compact set, the uniform norm is denoted by $\left\Vert \cdot \right\Vert_\infty$. The set of Lebesgue squared integrable functions from $(0,l)$ to $\mathbb{R}$ is denoted by $L^2(0,l)$ and is a Hilbert space when endowed with its natural inner product $\langle f,g \rangle_{L^2(0,l)} = \int_0^l f(y)g(y) \mathrm{d}y$. The associated norm is denoted by $\norme{\cdot}_{L^2(0,l)}$. For any $m\in\mathbb{N}$, $H^m(0,l)$ denotes the usual Sobolev space. Denoting by $\mathrm{AC}[0,l]$ the set of all absolutely continuous functions on $[0,l]$, $H^1(0,l) \subset \mathrm{AC}[0,l]$ in the sense that for any $f\in H^1(0,l)$, there exists a unique absolutely continuous function $g\in\mathrm{AC}[0,l]$ such that $f=g$ almost everywhere (in the sense of the Lebesgue measure), implying $f=g$ in $H^1(0,l)$. For a given normed vector space $(E,\norme{\cdot}_{E})$, $\mathcal{L}(E)$ denotes the space of bounded linear transformations from $E$ to $E$ and is a normed space when equipped with the induced norm denoted by $\left|\left|\cdot\right|\right|$. The range of a given operator $\mathcal{A}$ is denoted by $R(\mathcal{A})$ while its resolvent set is denoted by $\rho(\mathcal{A})$ and its kernel is defined by $ \mathrm{ker}(\mathcal{A})=\mathcal{A}^{-1}(\{0\})$. Further details can be found in, e.g., \cite[Annex A]{Curtain2012} and~\cite{Leoni2009}.

When dealing with the abstract form, the time derivative of a real-valued differentiable function $f:\mathbb{R}_+\rightarrow\mathbb{R}$ is denoted by $\dot{f}$. If $\mathcal{H}$ is a Hilbert space, the time derivative of a $\mathcal{H}$-valued differentiable function $f:\mathbb{R}_+\rightarrow\mathcal{H}$ is denoted by $\mathrm{d} f / \mathrm{d}t$.

\section{Problem Formulation and Boundary Control Law}\label{sec: model}

\subsection{Flexible wing model}
Let $l\in\mathbb{R}_+^*$ be the length of the wing. The structural parameters of the wing, assumed to be functions of the spatial variable $y$, are the mass per unit of span $\rho\in L^\infty(0,l)$, the moment of inertia per unit length $I_w\in L^\infty(0,l)$, and the bending (resp. torsional) stiffness $EI\in L^\infty(0,l)$ (resp. $GJ\in L^\infty(0,l)$). The damping characteristics of the wing are represented by the bending (resp. torsional) Kelvin-Voigt damping coefficient $\eta_\omega\in L^\infty(0,l)$ (resp. $\eta_\phi\in L^\infty(0,l)$). It is assumed that the essential infimum of these parameters over the wingspan are strictly positive, i.e., $\underline{\rho},\underline{I_w},\underline{EI},\underline{GJ},\underline{\eta_\omega},\underline{\eta_\phi} > 0$. 

To describe the dynamics of the flexible wing, we introduce $\omega : [0,l] \times \mathbb{R}_+ \rightarrow \mathbb{R}$ and $\phi : [0,l] \times \mathbb{R}_+ \rightarrow \mathbb{R}$ which denote, respectively, the bending and twisting displacements of the wing along the wingspan. The dynamics of the flexible wing are described by the following set of PDEs, which is a linear version of~\cite{Bialy2016} obtained by neglecting the elastic axis offset:
\begin{subequations}
\begin{align}
\rho \omega_{tt} + \left( EI \omega_{yy} + \eta_\omega EI \omega_{tyy} \right)_{yy} & = \rho \left( \alpha_\omega \phi + \beta_\omega \phi_t + \gamma_\omega \omega_t \right) , \label{eq: EDP flexion}  \\ 
I_w \phi_{tt} - \left( GJ \phi_{y} + \eta_\phi GJ \phi_{ty} \right)_{y} & = I_w \left( \alpha_\phi \phi + \beta_\phi \phi_t + \gamma_\phi \omega_t \right) , \label{eq: EDP torsion} 
\end{align}
\end{subequations}
in $(0,l)\times\mathbb{R}_+$, where $\alpha_\omega,\beta_\omega,\gamma_\omega,\alpha_\phi,\beta_\phi,\gamma_\phi\in L^\infty(0,l)$ represent the aerodynamic coefficients which are functions of the spatial variable $y$. The boundary conditions are such that, for any $t\geq0$, 
\begin{subequations}
\begin{align}
\omega(0,t)=\omega_{y}(0,t)=\phi(0,t) = & 0 , \label{eq: EDP boundary zero 1} \\
(EI \omega_{yy} + \eta_\omega EI \omega_{tyy}) (l,t) = & 0 , \label{eq: EDP boundary zero 2} \\
( EI \omega_{yy} + \eta_\omega EI \omega_{tyy} )_{y}(l,t) = & - L_\mathrm{tip}(t) + m_s \omega_{tt}(l,t) , \label{eq: EDP boundary flexion} \\ 
( GJ \phi_{y} + \eta_\phi GJ \phi_{ty} )(l,t) = & M_\mathrm{tip}(t) - J_s \phi_{tt}(l,t) , \label{eq: EDP boundary torsion}
\end{align} 
\end{subequations}
where $L_\mathrm{tip}:\mathbb{R}_+\rightarrow\mathbb{R}$ and $M_\mathrm{tip}:\mathbb{R}_+\rightarrow\mathbb{R}$ denote the control inputs located at the wing tip. Physically, $L_\mathrm{tip}(t)$ and $M_\mathrm{tip}(t)$ represent the aerodynamic lift force and pitching moment generated at time $t$ by the flaps located at the wing tip. The store at the wing tip is characterized by its mass $m_s\in\mathbb{R}_+^*$ and its moment of inertia $J_s\in\mathbb{R}_+^*$.  

Finally, the initial conditions are assumed to be: $\omega(\cdot,0) = \omega_0$ , $\omega_t(\cdot,0) = \omega_{t0}$, $\phi(\cdot,0) = \phi_0$, $\phi_t(\cdot,0) = \phi_{t0}$.

\subsection{Boundary control law}
For control design and practical implementation purposes, we make the following assumption. 

\begin{assum}
It is assumed that $\omega(l,\cdot)$, $\omega_{t}(l,\cdot)$, $\omega_{tt}(l,\cdot)$, $\phi(l,\cdot)$, $\phi_{t}(l,\cdot)$, and $\phi_{tt}(l,\cdot)$ are measured at the wing tip and available for feedback control.
\end{assum}

The proposed boundary stabilization control takes the following form: 
\begin{subequations}
\begin{align}
L_\mathrm{tip}(t) & = - k_1 \left[ \omega_{t}(l,t) + \epsilon_1 \omega(l,t) \right] + m_s \omega_{tt}(l,t) + u_1(t) , \label{eq: bf L_tip} \\
M_\mathrm{tip}(t) & = - k_2 \left[\phi_{t}(l,t) + \epsilon_2 \phi(l,t) \right] + J_s \phi_{tt}(l,t) + u_2(t) , \label{eq: bf M_tip}
\end{align}
\end{subequations}
for any $t\geq0$, where $k_1,k_2 \in\mathbb{R}_+$ are tunable controller gains that can be freely selected while $\epsilon_1,\epsilon_2 \in\mathbb{R}_+^*$ are two parameters that will be determined later in order to ensure adequate properties for the closed-loop system. The signals $u_1,u_2\in\mathcal{C}^2(\mathbb{R}_+,\mathbb{R})$ can be either auxiliary control inputs or disturbance inputs. In the remainder of this note, we study the stability properties of the system in closed loop with the proposed boundary control strategy.

\subsection{Closed-loop system in abstract form}
To analyze the properties of the closed-loop system with the proposed boundary control law, the problem is rewritten in abstract form. First, the following real Hilbert space is introduced:
\begin{equation*}
\begin{split}
\mathcal{H} = \{ & (f,g,h,z)\in H^2(0,l) \times L^2(0,l) \times H^1(0,l) \times L^2(0,l) \, : \\
& f(0)=f'(0)=0 , \, h(0)=0 \},
\end{split}
\end{equation*}
endowed with the inner product $\left\langle \cdot,\cdot \right\rangle_{\mathcal{H},1}$ defined by
\begin{equation*}
\begin{split}
& \langle (f_1,g_1,h_1,z_1),(f_2,g_2,h_2,z_2) \rangle_{\mathcal{H},1} \\
& \phantom{====} \triangleq \int_0^l [ EI(y) f_1''(y) f_2''(y) + \rho(y) g_1(y) g_2(y) ] \mathrm{d}y \\
& \phantom{==== \triangleq} + \int_0^l [ GJ(y) h_1'(y) h_2'(y) + I_w(y) z_1(y) z_2(y) ] \mathrm{d}y .
\end{split}
\end{equation*}
Note that $\left\langle \cdot,\cdot \right\rangle_{\mathcal{H},1}$ is indeed an inner product on $\mathcal{H}$ because $\underline{EI},\underline{GJ},\underline{\rho},\underline{I_w}>0$. The motivation for introducing this specific inner product is closely related to the physical nature of the wing energy (\ref{eq: definition energy}). Indeed, denoting by  $\norme{\cdot}_{\mathcal{H},1}$ the induced norm, the energy of the wing due to elastic deformations, including kinetic and potential energy, is given by 
\begin{equation}\label{eq: definition energy}
\forall t\geq 0, \;E(t) = \dfrac{1}{2} \norme{(\omega(\cdot,t) , \omega_t(\cdot,t) , \phi(\cdot,t) , \phi_t(\cdot,t))}_{\mathcal{H},1}^2 .
\end{equation}

To study the well-posedness of the closed-loop system in the presence of input perturbations, we define the following abstract operator:
\begin{equation}\label{eq: def U}
\begin{array}{lccc}
\mathcal{A}_d : & D(\mathcal{A}_d) & \longrightarrow & \mathcal{H} \\
& (f,g,h,z) & \longrightarrow & (g,\tilde{g},z,\tilde{z})
\end{array}
\end{equation}
where
\begin{align*}
\tilde{g} & \triangleq - \dfrac{1}{\rho} (EI f''+\eta_\omega EI g'')'' + \alpha_\omega h + \beta_\omega z + \gamma_\omega g , \\
\tilde{z} & \triangleq \dfrac{1}{I_w} (GJ h' + \eta_\phi GJ z')' + \alpha_\phi h + \beta_\phi z + \gamma_\phi g ,
\end{align*}
with domain $D(\mathcal{A}_d) \subset \mathcal{H}$ defined by
\begin{equation}\label{eq: def domain U}
\begin{split}
D(\mathcal{A}_d) \triangleq \{ & (f,g,h,z)\in\mathcal{H} \, : \\
& \, g \in H^2(0,l) ,\, z \in H^1(0,l) ,\, \\
& EI f'' + \eta_\omega EI g'' \in H^2(0,l) ,\, \\
& GJ h' + \eta_\phi GJ z' \in H^1(0,l) ,\,  \\
& f(0)=f'(0)=0 ,\, g(0)=g'(0)=0 ,\, \\
& h(0)=0 ,\, z(0)=0 ,\, \\
& (EI f'' + \eta_\omega EI g'')(l)=0 \} \\
\end{split}
\end{equation}
We also introduce the boundary operator:
\begin{equation}
\begin{array}{lccc}
\mathcal{B} : & D(\mathcal{B}) & \longrightarrow & \mathbb{R}^2 \\
& (f,g,h,z) & \longrightarrow & (\tilde{u}_1,\tilde{u}_2)
\end{array}
\end{equation}
where $\mathbb{R}^2$ is endowed with the usual 2-norm,
\begin{align*}
\tilde{u}_1 & \triangleq - ( EI f'' + \eta_\omega EI g'')'(l) + k_1 (g(l)+\epsilon_1 f(l)) , \\
\tilde{u}_2 & \triangleq (GJ h' + \eta_\phi GJ z')(l) + k_2(z(l) + \epsilon_2 h(l)) ,
\end{align*}
with domain $D(\mathcal{B})\triangleq D(\mathcal{A}_d)\subset\mathcal{H}$. Let  $U = (u_1,u_2) \in \mathcal{C}^2(\mathbb{R}_+,\mathbb{R}^2)$ be the disturbing input. It leads to the following abstract boundary control problem:
\begin{equation}\label{eq: boundary control problem}
\left\{\begin{split}
\dfrac{\mathrm{d} X}{\mathrm{d} t}(t) & = \mathcal{A}_d X(t) ,\; t>0 \\
\mathcal{B} X(t) & = U(t) ,\; t \geq 0 \\
X(0) & = X_0 \in D(\mathcal{A}_d) \;\;\mathrm{s.t.}\;\; \mathcal{B} X_0 = U(0)
\end{split}\right.
\end{equation}
where $X(t) = \left( \omega(\cdot,t) , \omega_t(\cdot,t) , \phi(\cdot,t) , \phi_t(\cdot,t) \right)$ is the state vector and $X_0 = \left( \omega_0 , \omega_{t0} , \phi_0 , \phi_{t0} \right)$ is the initial condition. To study the stability properties of the boundary control problem (\ref{eq: boundary control problem}), its well-posedness is first investigated in the next section.

\section{Well-posedness assessment}\label{sec: well-posedness}
In order to study the well-posedness of the boundary control problem (\ref{eq: boundary control problem}), it is useful to first study the disturbance free version of (\ref{eq: boundary control problem}), i.e., for $U=0$. To do so, we introduce the associated operator $\mathcal{A} \triangleq \left. \mathcal{A}_d \right|_{D(\mathcal{A})}$ with $D(\mathcal{A}) \triangleq  D(\mathcal{A}_d) \cap \mathrm{ker}(\mathcal{B})$. To facilitate the upcoming developments, the following two linear operators $\mathcal{A}_1 : D(\mathcal{A}_1) \rightarrow \mathcal{H}$ and $\mathcal{A}_2 : D(\mathcal{A}_2) \rightarrow \mathcal{H}$ are introduced:
\begin{align*}
\mathcal{A}_1 & (f,g,h,z) \\
\triangleq & \left( g , -\dfrac{1}{\rho}(EI f'' + \eta_\omega EI g'')'' , z , \dfrac{1}{I_w} (GJ h' + \eta_\phi GJ z')' \right) , \\
\mathcal{A}_2 & (f,g,h,z)  \triangleq \left( 0 , \alpha_\omega h + \beta_\omega z + \gamma_\omega g , 0 , \alpha_\phi h + \beta_\phi z + \gamma_\phi g \right) ,
\end{align*}
with domains $D(\mathcal{A}_1)=D(\mathcal{A})$ and $D(\mathcal{A}_2)=\mathcal{H}$. Obviously, $\mathcal{A}=\mathcal{A}_1+\mathcal{A}_2$ over $D(\mathcal{A})$.

The following two inequalities will be used in the subsequent developments.
\begin{lemma}\label{lemma: Poincare's and Agmon's inequalities}~\cite{Hardy1952,Krstic2008}
For any $f\in H^1(0,l)\subset \mathrm{AC}(0,l)$\footnote{Inclusion in the sense explained in the introduction.} such that $f(0)=0$, the Poincar\'{e}'s inequality ensures that
\begin{equation*}
\norme{f}_{L^2(0,l)}^2 \leq \dfrac{4 l^2}{\pi^2} \norme{f'}_{L^2(0,l)}^2 ,
\end{equation*}
while the Agmon's inequality provides
\begin{equation*}
\norme{f}_\infty^2 \leq 2 \norme{f}_{L^2(0,l)} \norme{f'}_{L^2(0,l)} .
\end{equation*}
\end{lemma}

\subsection{Necessity and introduction of a second inner product on $\mathcal{H}$}
The following Lemma shows that $\mathcal{A}_1$ is not dissipative with respect to the inner product $\left\langle \cdot,\cdot \right\rangle_{\mathcal{H},1}$ and hence, the Lumer-Philips theorem~\cite{Curtain2012,Pazy2012} implies that $\mathcal{A}_1$ does not generate a $C_0$-semigroup of contractions on $(\mathcal{H},\left\langle \cdot,\cdot \right\rangle_{\mathcal{H},1})$. 

\begin{lemma}
The operator $\mathcal{A}_1$ is not dissipative on $\mathcal{H}$ endowed with $\left\langle \cdot,\cdot \right\rangle_{\mathcal{H},1}$.
\end{lemma}

\textbf{Proof.}
Integrating by parts, we have for any $X = (f,g,h,z)\in D(\mathcal{A}_1)$, 
\begin{align}
& \left\langle \mathcal{A}_1 X , X \right\rangle_{\mathcal{H},1} \nonumber \\
= & - k_1 (g(l)+\epsilon_1 f(l)) g(l) - \int_0^l \eta_\omega(y) EI(y) g''(y)^2 \mathrm{d}y  \label{eq: dissipativity A1 2} \\
& - k_2 (z(l)+\epsilon_2 h(l))z(l) - \int_0^l \eta_\phi(y) GJ(y) z'(y)^2 \nonumber \mathrm{d}y \nonumber .
\end{align}
In particular, considering $f=g=0$ and, for all $y \in [0,l]$, 
\begin{equation*}
h(y) = \int_0^y \dfrac{\kappa_1 \xi + \kappa_2}{GJ(\xi)}\mathrm{d}\xi , \;\; 
z(y)=\int_0^y\dfrac{\kappa_3}{\eta_\phi(\xi) GJ(\xi)} \mathrm{d}\xi ,
\end{equation*}
where\footnote{Note that $\kappa_1$ and $\kappa_2$ are well defined because $1/GJ(y) \geq 1/\overline{GJ}$ for almost all $y \in [0,l]$, which implies that $l I_1 - I_2 \geq l^2/(2 \overline{GJ})>0$.}
\begin{equation*}
\kappa_1 = \dfrac{1}{l I_1 - I_2} \left\{ \dfrac{1}{k_2 \epsilon_2} \left( 1 + k_2 + \dfrac{1}{I_3} \right) + I_1 \right\} , 
\end{equation*}
\begin{equation*}
\kappa_2 = \dfrac{-1}{l I_1 - I_2} \left\{ \dfrac{l}{k_2 \epsilon_2} \left( 1 + k_2 + \dfrac{1}{I_3} \right) + I_2 \right\} , \;\;
\kappa_3 = \dfrac{1}{I_3} ,
\end{equation*}
with
\begin{equation*}
I_1 = \int_0^l \dfrac{\mathrm{d}\xi}{GJ(\xi)} ,\; I_2 = \int_0^l \dfrac{\xi}{GJ(\xi)} \mathrm{d}\xi ,\; I_3 = \int_0^l \dfrac{\mathrm{d}\xi}{\eta_\phi(\xi) GJ(\xi)}  ,
\end{equation*}
we have $X = (f,g,h,z)\in D(\mathcal{A})$ and, based on (\ref{eq: dissipativity A1 2}), straightforward calculations yields $\left\langle \mathcal{A}_1 X , X \right\rangle_{\mathcal{H},1} = 1 > 0$. Hence, $\mathcal{A}_1$ is not dissipative relatively to $\left\langle \cdot,\cdot \right\rangle_{\mathcal{H},1}$.
\qed

To solve this problem, there exist at least two possible alternatives. The first one is to resort to the Hille-Yosida theorem~\cite{Curtain2012,Pazy2012} to ensure that $\mathcal{A}_1$ generates a $C_0$-semigroup on $(\mathcal{H},\left\langle \cdot,\cdot \right\rangle_{\mathcal{H},1})$, which is a weaker property than the $C_0$-semigroup of contractions. The second one, adopted in this work, is still to apply the Lumer-Philips theorem while considering another inner product $\left\langle \cdot,\cdot \right\rangle_{\mathcal{H},2}$ on $\mathcal{H}$.

Let $\epsilon_1,\epsilon_2\in\mathbb{R}_+^*$ be constant parameters with constraints given later in Lemma~\ref{lemma: inner product 2} and $\left\langle\cdot,\cdot\right\rangle_{\mathcal{H},2}:\mathcal{H}\times\mathcal{H} \rightarrow \mathbb{R}$ be defined for any $(f_1,g_1,h_1,z_1),(f_2,g_2,h_2,z_2)\in\mathcal{H}$ by
\begin{equation*}
\begin{split}
& \langle (f_1,g_1,h_1,z_1) ,(f_2,g_2,h_2,z_2) \rangle_{\mathcal{H},2} \\
& \phantom{=====} \triangleq \left\langle (f_1,g_1,h_1,z_1),(f_2,g_2,h_2,z_2) \right\rangle_{\mathcal{H},1} \\
& \phantom{=====\triangleq} + \epsilon_1 \int_0^l \rho(y) \left[ f_1(y) g_2(y) + g_1(y) f_2(y) \right] \mathrm{d}y \\
& \phantom{=====\triangleq} + \epsilon_2 \int_0^l I_w(y) \left[ h_1(y) z_2(y) + z_1(y) h_2(y) \right] \mathrm{d}y .
\end{split}
\end{equation*}
We also introduce a constant $K_m\in\mathbb{R}_+^*$ defined by
\begin{equation*}
K_m = \max\left( \sqrt{\overline{\rho}} , \dfrac{16 l^4 \sqrt{\overline{\rho}}}{\pi^4 \underline{EI}} , \sqrt{\overline{I_w}} , \dfrac{4 l^2 \sqrt{\overline{I_w}}}{\pi^2 \underline{GJ}} \right) .
\end{equation*}
Then, the following lemma holds.

\begin{lemma}\label{lemma: inner product 2}
For any given $0<\epsilon_1,\epsilon_2<1/K_m$, $\left\langle\cdot,\cdot\right\rangle_{\mathcal{H},2}$ is an inner product for $\mathcal{H}$. Furthermore, the norm induced from this inner product, denoted by $\norme{\cdot}_{\mathcal{H},2}$, is equivalent to $\norme{\cdot}_{\mathcal{H},1}$. Thus, $(\mathcal{H},\left\langle\cdot,\cdot\right\rangle_{\mathcal{H},2})$ is a real Hilbert space.
\end{lemma}

\textbf{Proof.} First, $\left\langle\cdot,\cdot\right\rangle_{\mathcal{H},2}$ is bilinear and symmetric. For any $X=(f,g,h,z)\in\mathcal{H}$, 
\begin{align*}
\left\langle X , X \right\rangle_{\mathcal{H},2}
= & \norme{X}_{\mathcal{H},1}^2
+ 2\epsilon_1 \int_0^l \rho(y) f(y) g(y) \mathrm{d}y \\
& + 2\epsilon_2 \int_0^l I_w(y) h(y) z(y) \mathrm{d}y .
\end{align*}
Then, by first applying Young's inequality\footnote{For any $a,b\in\mathbb{R}_+$ and $r\in\mathbb{R}_+^*$, the Young's inequality provides $ab \leq a^2/(2r) + r b^2/2$.}, and then Poincar\'{e}'s inequality, one has for any $X\in\mathcal{H}$, 
\begin{equation}\label{eq: norm equivalent}
( 1 - \epsilon_m K_m ) \norme{X}_{\mathcal{H},1}^2
\leq 
\left\langle X , X \right\rangle_{\mathcal{H},2} 
\leq 
( 1 + \epsilon_m K_m ) \norme{X}_{\mathcal{H},1}^2 ,
\end{equation}
with $\epsilon_m = \max(\epsilon_1,\epsilon_2)$. Then, for $0<\epsilon_m<1/K_m$, $\left\langle\cdot,\cdot\right\rangle_{\mathcal{H},2}$ is positive and definite and hence, it defines an inner product for $\mathcal{H}$. Furthermore, denoting by $\norme{\cdot}_{\mathcal{H},2}$ the associated norm, (\ref{eq: norm equivalent}) implies that $\norme{\cdot}_{\mathcal{H},2}$ and $\norme{\cdot}_{\mathcal{H},1}$ are equivalent. It follows that $(\mathcal{H},\left\langle\cdot,\cdot\right\rangle_{\mathcal{H},2})$ is a real Hilbert space.
\qed

In the subsequent developments, we assume that the controller parameters $\epsilon_1$ and $\epsilon_2$ are constrained by $0<\epsilon_1,\epsilon_2<1/K_m$. Therefore, Lemma~\ref{lemma: inner product 2} is applied hereafter.

\subsection{$\mathcal{A}_1$ generates a $C_0$-semigroup of contractions on $(\mathcal{H},\left\langle \cdot,\cdot \right\rangle_{\mathcal{H},2})$}
To apply the Lumer-Phillips theorem in the context of an Hilbert space, we need to assess a dissipativity condition and a certain range condition~\cite{Curtain2012,Pazy2012}. To assess the first point, we introduce the following two constants:
\begin{align*}
\epsilon_1^* = \dfrac{4 \pi^4 \underline{\eta_\omega EI}}{64 l^4 \overline{\rho} + \pi^4 \overline{\eta_\omega} \underline{\eta_\omega EI}} ,\; \epsilon_2^* = \dfrac{4 \pi^2 \underline{\eta_\phi GJ}}{16 l^2 \overline{I_w} + \pi^2 \overline{\eta_\phi} \underline{\eta_\phi GJ}} . 
\end{align*}

\begin{lemma} \label{lemma: dissipative A1}
Let $\epsilon_1,\epsilon_2 \in \mathbb{R}_+^*$ such that $\epsilon_1 < \min( \epsilon_1^* , 1/K_m)$ and $\epsilon_2 < \min( \epsilon_2^* , 1/K_m)$. Then, the operator $\mathcal{A}_1 : D(\mathcal{A}_1) \rightarrow \mathcal{H}$ is dissipative with respect to $\left\langle  \cdot , \cdot \right\rangle_{\mathcal{H},2}$.
\end{lemma}

\textbf{Proof.} As $0<\epsilon_1,\epsilon_2<1/K_m$, Lemma~\ref{lemma: inner product 2} ensures that $(\mathcal{H},\left\langle  \cdot , \cdot \right\rangle_{\mathcal{H},2})$ is a real Hilbert space. Letting $X = (f,g,h,z) \in D(\mathcal{A}_1)$, based on (\ref{eq: dissipativity A1 2}) and integrations by parts, it yields,
\begin{align*}
& \left\langle \mathcal{A}_1 X , X \right\rangle_{\mathcal{H},2} \\
= & - k_1 (g(l)+\epsilon_1 f(l))^2 - k_2 (z(l)+\epsilon_2 h(l))^2 \\
& + \epsilon_1 \int_0^l \rho(y) g(y)^2 \mathrm{d}y + \epsilon_2 \int_0^l I_w(y) z(y)^2 \mathrm{d}y \\
& - \int_0^l \eta_\omega(y) EI(y) g''(y)^2 \mathrm{d}y - \int_0^l \eta_\phi(y) GJ(y) z'(y)^2 \mathrm{d}y \\
& - \epsilon_1 \int_0^l EI(y) f''(y)^2 \mathrm{d}y - \epsilon_1 \int_0^l \eta_\omega(y) EI(y) f''(y)  g''(y) \mathrm{d}y \\
& - \epsilon_2 \int_0^l GJ(y) h'(y)^2 \mathrm{d}y - \epsilon_2 \int_0^l \eta_\phi(y) GJ(y) h'(y)  z'(y) \mathrm{d}y 
\end{align*}
First applying Young's inequality and then Poincar{\'e}'s inequality, it provides for all $X=(f,g,h,z)\in D(\mathcal{A})$ and for all $r_1,r_2>0$,
\begin{align}
\left\langle \mathcal{A}_1 X , X \right\rangle_{\mathcal{H},2} \leq & - k_1 (g(l)+\epsilon_1 f(l))^2 - k_2 (z(l)+\epsilon_2 h(l))^2 \nonumber \\
& - \left( 1 - \dfrac{\epsilon_1}{\varphi_1(r_1)} \right) \int_0^l \eta_\omega(y) EI(y) g''(y)^2 \mathrm{d}y \nonumber \\
& - \left( 1 - \dfrac{\epsilon_2}{\varphi_2(r_2)} \right) \int_0^l \eta_\phi(y) GJ(y) z'(y)^2 \mathrm{d}y \nonumber \\
& - \epsilon_1 \left( 1 - \dfrac{\sqrt{\overline{\eta_\omega}}}{2 r_1} \right) \int_0^l EI(y) f''(y)^2 \mathrm{d}y \nonumber \\
& - \epsilon_2 \left( 1 - \dfrac{\sqrt{\overline{\eta_\phi}}}{2 r_2} \right) \int_0^l GJ(y) h'(y)^2 \mathrm{d}y \label{eq: inner A1X-X}
\end{align}
where $\varphi_1: \mathbb{R}_+^* \ni x \rightarrow 2\pi^4\underline{\eta_\omega EI}/(32 l^4 \overline{\rho} +\pi^4\sqrt{\overline{\eta_\omega}}\underline{\eta_\omega EI} x)$ and $\varphi_2: \mathbb{R}_+^* \ni x \rightarrow 2 \pi^2 \underline{\eta_\phi GJ} / (8 l^2 \overline{I_w} + \pi^2 \sqrt{\overline{\eta_\phi}} \underline{\eta_\phi GJ} x)$. As $\varphi_1$ is a continuous decreasing function over $\mathbb{R}_+^*$ and, by assumption, $\epsilon_1 < \epsilon_1^* = \varphi_1(\sqrt{\overline{\eta_\omega}}/2)$, there exists $r_1^*>\sqrt{\overline{\eta_\omega}}/2$ such that $\epsilon_1 < \varphi_1(r_1^*) < \varphi_1(\sqrt{\overline{\eta_\omega}}/2)$. Similarly, there exists $r_2^*>\sqrt{\overline{\eta_\phi}}/2$ such that $\epsilon_2 < \varphi_2(r_2^*) < \varphi_2(\sqrt{\overline{\eta_\phi}}/2) = \epsilon_2^*$. Therefore, taking $r_1=r_1^*$ and $r_2=r_2^*$ in (\ref{eq: inner A1X-X}), it ensures that for all $X\in\mathcal{H}$, $\left\langle \mathcal{A}_1 X , X \right\rangle_{\mathcal{H},2} \leq 0$, i.e., $\mathcal{A}_1$ is dissipative on $\mathcal{H}$ endowed with $\left\langle  \cdot , \cdot \right\rangle_{\mathcal{H},2}$.
\qed

We now investigate the range condition.

\begin{lemma} \label{lemma: inv A1}
The operator $\mathcal{A}_1^{-1} : \mathcal{H} \rightarrow D(\mathcal{A}_1)$ exists and is bounded, i.e., $\mathcal{A}_1^{-1}\in\mathcal{L}(\mathcal{H})$. Hence, $0\in\rho(\mathcal{A}_1)$ and $\mathcal{A}_1$ is a closed operator.
\end{lemma}

\textbf{Proof.} 
We first investigate the surjectivity of $\mathcal{A}_1$. Let $(\tilde{f},\tilde{g},\tilde{h},\tilde{z})\in\mathcal{H}$. Then, with $f$ defined by
\begin{align}
f(y) = & - \int_0^y ( y - \xi ) \eta_\omega(\xi) \tilde{f}''(\xi) \mathrm{d}\xi \label{eq: surjectivite operator A1 def f} \\
&  - k_1 \int_0^y \dfrac{(y-\xi)(l-\xi)}{EI(\xi)} \mathrm{d}\xi \, (\tilde{f}(l)+\epsilon_1 \alpha(\tilde{f},\tilde{g})) \nonumber \\
&  - \int_0^y \dfrac{y-\xi_1}{EI(\xi_1)} \int_{\xi_1}^l (\xi_2-\xi_1) \rho(\xi_2) \tilde{g}(\xi_2) \mathrm{d}\xi_2 \mathrm{d}\xi_1 \nonumber ,
\end{align}
where
\begin{align*}
& \alpha(\tilde{f},\tilde{g}) \\ 
= & - \left\{ 1 + k_1 \epsilon_1 \int_0^{l} \dfrac{(l-\xi)^2}{EI(\xi)} \mathrm{d}\xi \right\}^{-1} \\
& \times 
\left\{ \int_0^l (l-\xi) \eta_\omega(\xi) \tilde{f}''(\xi) \mathrm{d}\xi + k_1 \int_0^l \dfrac{(l-\xi)^2}{EI(\xi)} \mathrm{d}\xi \, \tilde{f}(l) \right. \\
& \phantom{\times \{\,\;\;\; } + \left. \int_0^l \dfrac{l-\xi_1}{EI(\xi_1)} \int_{\xi_1}^l (\xi_2-\xi_1) \rho(\xi_2) \tilde{g}(\xi_2) \mathrm{d}\xi_2 \mathrm{d}\xi_1 \right\} ,
\end{align*}
$g=\tilde{f}$, $h$ defined by
\begin{align}
& h(y) \label{eq: surjectivite operator A1 def h} \\
= & - \int_0^y \eta_\phi(\xi) \tilde{h}'(\xi) \mathrm{d}\xi - k_2 \int_0^y \dfrac{\mathrm{d}\xi}{GJ(\xi)} \, (\tilde{h}(l)+\epsilon_2 \beta(\tilde{h},\tilde{z})) \nonumber \\
&  - \int_0^y \dfrac{1}{GJ(\xi_1)} \int_{\xi_1}^l I_w(\xi_2) \tilde{z}(\xi_2) \mathrm{d}\xi_2 \mathrm{d}\xi_1 \nonumber ,
\end{align}
where
\begin{align*}
\beta(\tilde{h},\tilde{z})
& = - \left\{ 1 + k_2 \epsilon_2 \int_0^{l} \dfrac{\mathrm{d}\xi}{GJ(\xi)} \mathrm{d}\xi \right\}^{-1} \\
& \times 
\left\{ \int_0^l \eta_\phi(\xi) \tilde{h}'(\xi) \mathrm{d}\xi + k_2 \int_0^l \dfrac{\mathrm{d}\xi}{GJ(\xi)} \, \tilde{h}(l) \right. \nonumber \\
& \phantom{\times \{\,\;\;\; } \left. + \int_0^l \dfrac{1}{GJ(\xi_1)} \int_{\xi_1}^l I_w(\xi_2) \tilde{z}(\xi_2) \mathrm{d}\xi_2 \mathrm{d}\xi_1 \right\} ,
\end{align*}
and $z=\tilde{h}$, we have $(f,g,h,z)\in D(\mathcal{A}_1)=D(\mathcal{A})$ and $(\tilde{f},\tilde{g},\tilde{h},\tilde{z}) = \mathcal{A}_1(f,g,h,z)$. Thus, $\mathcal{A}_1$ is onto $\mathcal{H}$. In addition, as $(f,g,h,z)$ depends linearly on $(\tilde{f},\tilde{g},\tilde{h},\tilde{z})\in\mathcal{H}$, it shows that the operator $\mathcal{A}_1$ is right invertible. To conclude that $\mathcal{A}_1$ is invertible, we investigate the injectivity. Let $(f,g,h,z)\in D(\mathcal{A}_1)$ such that $\mathcal{A}_1(f,g,h,z)=(0,0,0,0)$. Thus, $g=z=0$, which yields $(EI f'')''=0$ and $(GJ h')'=0$ along with $f(0)=f'(0)=(EI f'')(l)=0$, $(EI f'')'(l) = k_1 \epsilon_1 f(l)$, $h(0)=0$, and $(GJ h')(l) = -k_2 \epsilon_2 h(l)$. Then, as $(EI f'')',(EI f''), f',f \in \mathrm{AC}[0,l]$, it yields for any $y\in[0,l]$ by successive integrations,
\begin{equation*}
f(y) = - k_1 \epsilon_1 \int_0^y \dfrac{(y-\xi)(l-\xi)}{EI(\xi)} \mathrm{d}\xi \, f(l) .
\end{equation*}
Evaluating at $y=l$,
\begin{equation*}
\underbrace{\left( 1 + k_1 \epsilon_1 \int_0^l \dfrac{(l-\xi)^2}{EI(\xi)} \mathrm{d}\xi \right)}_{>0} f(l) = 0
\Rightarrow
f(l) = 0 ,
\end{equation*}
which implies that $f=0$. Similarly, one can show that $h=0$. Hence, the operator $\mathcal{A}_1$ is injective.

Therefore, $\mathcal{A}_1^{-1} : \mathcal{H} \rightarrow D(\mathcal{A}_1)$ is well defined for any $(\tilde{f},\tilde{g},\tilde{h},\tilde{z})\in\mathcal{H}$ by $\mathcal{A}_1^{-1}(\tilde{f},\tilde{g},\tilde{h},\tilde{z}) = (f,\tilde{f},h,\tilde{h})$ where $f$ and $h$ are given by (\ref{eq: surjectivite operator A1 def f}) and (\ref{eq: surjectivite operator A1 def h}), respectively. Based on Poincar{\'e}'s inequality, straightforward computations show that $\mathcal{A}_1^{-1}$ is a bounded operator. Hence, it ensures that $\mathcal{A}_1^{-1}\in\mathcal{L}(\mathcal{H})$, i.e., $0\in\rho(\mathcal{A}_1)$ and $\mathcal{A}_1$ is a closed operator.
\qed

We can now introduce the following property regarding $\mathcal{A}_1$.

\begin{theorem}\label{th: A1 generation semigroup}
Let $\epsilon_1,\epsilon_2 \in \mathbb{R}_+^*$ such that $\epsilon_1 < \min( \epsilon_1^* , 1/K_m)$ and $\epsilon_2 < \min( \epsilon_2^* , 1/K_m)$. Then, the operator $\mathcal{A}_1$ generates a $C_0$-semigroup of contractions on $(\mathcal{H},\left\langle  \cdot , \cdot \right\rangle_{\mathcal{H},2})$. Furthermore, $D(\mathcal{A}_1)$ is dense in $\mathcal{H}$ endowed by either $\langle\cdot,\cdot\rangle_{\mathcal{H},1}$ or $\langle\cdot,\cdot\rangle_{\mathcal{H},2}$.
\end{theorem}

\textbf{Proof.} Under the assumptions of the theorem, $(\mathcal{H},\left\langle  \cdot , \cdot \right\rangle_{\mathcal{H},2})$ is a real Hilbert space and $\mathcal{A}_1$ is dissipative with respect to $\left\langle  \cdot , \cdot \right\rangle_{\mathcal{H},2}$. Furthermore, as the resolvent set of a closed operator is an open subset of $\mathbb{C}$, $0\in\rho(\mathcal{A}_1)$ implies the existence of $\lambda_0>0$ such that $\lambda_0 \in \rho(\mathcal{A}_1)$. In particular, $R(\lambda_0 I_{D(\mathcal{A}_1)} - \mathcal{A}_1) = \mathcal{H}$. Therefore, the application of the Lumer-Philips theorem for reflexive spaces~\cite[Th.2.29]{Luo2012}~\cite[Chap.1, Th.4.5 and Th.4.6]{Pazy2012} concludes the proof.
\qed

\subsection{$\mathcal{A}$ generates a $C_0$-semigroup}
The following lemma is a direct consequence of the definition of the operator $\mathcal{A}_2$ and the application of Young's and Poincar\'{e}'s inequalities.

\begin{lemma}\label{lemma: A2}
Operator $\mathcal{A}_2$ is bounded, i.e., $\mathcal{A}_2\in\mathcal{L}(\mathcal{H})$.
\end{lemma}

This result allows introducing the following main result.

\begin{theorem}\label{th: A generates C0-semigroup}
Let $\epsilon_1,\epsilon_2 \in \mathbb{R}_+^*$ such that $\epsilon_1 < \min( \epsilon_1^* , 1/K_m)$ and $\epsilon_2 < \min( \epsilon_2^* , 1/K_m)$. Then, the operator $\mathcal{A}$ generates a $C_0$-semigroup on $(\mathcal{H},\left\langle  \cdot , \cdot \right\rangle_{\mathcal{H},2})$.
\end{theorem}

\textbf{Proof.} Based on Theorem~\ref{th: A1 generation semigroup} and Lemma~\ref{lemma: A2}, the claimed result is a direct consequence of the perturbation theory~\cite[Chap.3, Th.1.1]{Pazy2012}\cite[Th.3.2.1]{Curtain2012}.
\qed

The following corollary is a consequence of the equivalence of the norms stated in Lemma~\ref{lemma: inner product 2} and the uniqueness of the $C_0$-semigroup associated to a given infinitesimal generator~\cite[Th.2.14]{Luo2012}.

\begin{corollary}\label{cor: C0 semi-group generation}
Under the assumptions of Theorem~\ref{th: A generates C0-semigroup}, $\mathcal{A}$ generates a $C_0$-semigroup on $(\mathcal{H},\left\langle  \cdot , \cdot \right\rangle_{\mathcal{H},1})$ which coincides with the $C_0$-semigroup generated by $\mathcal{A}$ on $(\mathcal{H},\left\langle  \cdot , \cdot \right\rangle_{\mathcal{H},2})$.
\end{corollary} 
 
\subsection{Well-posedness of the boundary control problem}
According to \cite[Def. 3.3.2]{Curtain2012}, we check that (\ref{eq: boundary control problem}) satisfies the definition of a \emph{boundary control system}. First, $D(\mathcal{A}) = D(\mathcal{A}_d) \cap \mathrm{ker}(\mathcal{B})$ and $\mathcal{A}z=\mathcal{A}_dz$ for any $z\in D(\mathcal{A})$. Second, we need to check the existence of an operator $B\in\mathcal{L}(\mathbb{R}^2,\mathcal{H})$ such that $R(B) \subset D(\mathcal{A}_d)$, $\mathcal{A}_dB\in\mathcal{L}(\mathbb{R}^2,\mathcal{H})$, and $\mathcal{B}B=I_{\mathbb{R}^2}$. We define the following candidate:
\begin{equation}
\begin{array}{lccc}
B : & \mathbb{R}^2 & \longrightarrow & \mathcal{H} \\
& (u_1,u_2) & \longrightarrow & (f_{u_1},0,h_{u_2},0)
\end{array}
\end{equation}
where functions $f_{u_1}$ and $h_{u_2}$ are defined for any $y\in[0,l]$ by
\begin{align*}
f_{u_1}(y) \triangleq & \dfrac{u_1}{b_1} \int_0^y \dfrac{(y-\xi)(l-\xi)}{EI(\xi)} \mathrm{d}\xi , \\
h_{u_2}(y) \triangleq & \dfrac{u_2}{b_2} \int_0^y \dfrac{\mathrm{d}\xi}{GJ(\xi)} ,
\end{align*}
with the constants:
\begin{equation*}
b_1 \triangleq 1 + k_1 \epsilon_1 \int_0^l \dfrac{(l - \xi)^2}{EI(\xi)} \mathrm{d}\xi  ,\;\; b_2 \triangleq 1+k_2 \epsilon_2 \int_0^l \dfrac{\mathrm{d}\xi}{GJ(\xi)} .
\end{equation*}
We check that all the required conditions are satisfied. First, $B$ is clearly linear and satisfies for any $U = (u_1,u_2) \in \mathbb{R}^2$,
\begin{align*}
\left\Vert B U \right\Vert_{\mathcal{H},1}^2
& = \dfrac{u_1^2}{b_1^2} \int_0^l \dfrac{(l-y)^2}{EI(y)} \mathrm{d}y + \dfrac{u_2^2}{b_2^2} \int_0^l \dfrac{\mathrm{d}y}{GJ(y)} \\
& \leq \max\left( \dfrac{1}{b_1^2} \int_0^l \dfrac{(l-y)^2}{EI(y)} \mathrm{d}y , \dfrac{1}{b_2^2} \int_0^l \dfrac{\mathrm{d}y}{GJ(y)} \right) \left\Vert U \right\Vert_2^2 ,
\end{align*}
where the equality holds for either $U = (1,0)$ or $U = (0,1)$. Thus $B\in\mathcal{L}(\mathbb{R}^2,\mathcal{H})$ with
\begin{equation}\label{eq: norm B}
\left\Vert B \right\Vert = \max\left( \dfrac{1}{b_1} \sqrt{\int_0^l \dfrac{(l-y)^2}{EI(y)} \mathrm{d}y} , \dfrac{1}{b_2} \sqrt{\int_0^l \dfrac{\mathrm{d}y}{GJ(y)}} \right) .
\end{equation} 
Furthermore, straightforward calculations show that for all $U\in\mathbb{R}^2$, $BU\in D(\mathcal{A}_d)$ and $\mathcal{B}BU=U$. Finally, for any $U = (u_1,u_2) \in \mathbb{R}^2$,
\begin{align*}
& \left\Vert \mathcal{A}_d B U \right\Vert_{\mathcal{H},1}^2 \\
= & \dfrac{u_2^2}{b_2^2} \int_0^l \left( \rho(y) \alpha_w(y)^2 + I_w(y) \alpha_\phi(y)^2 \right) \left( \int_0^y \dfrac{\mathrm{d}\xi}{GJ(\xi)} \right)^2 \mathrm{d}y \\
\leq & \dfrac{1}{b_2^2} \int_0^l \left( \rho(y) \alpha_w(y)^2 + I_w(y) \alpha_\phi(y)^2 \right) \left( \int_0^y \dfrac{\mathrm{d}\xi}{GJ(\xi)} \right)^2 \mathrm{d}y \left\Vert U \right\Vert_2^2 ,
\end{align*}
where the equality holds for $U = (0,1)$. Thus $\mathcal{A}_d B\in\mathcal{L}(\mathbb{R}^2,\mathcal{H})$ with
\begin{align}\label{eq: norm UB}
& \left\Vert \mathcal{A}_d B \right\Vert \\
= & \dfrac{1}{b_2} \sqrt{\int_0^l \left( \rho(y) \alpha_w(y)^2 + I_w(y) \alpha_\phi(y)^2 \right) \left( \int_0^y \dfrac{\mathrm{d}\xi}{GJ(\xi)} \right)^2 \mathrm{d}y} . \nonumber
\end{align} 

Based on Theorem~\ref{th: A generates C0-semigroup}, $\mathcal{A}$ generates a $C_0$-semigroup for sufficiently small parameters $\epsilon_1,\epsilon_2>0$. We deduce that for any $U \in \mathcal{C}^2(\mathbb{R}_+,\mathbb{R}^2)$, the abstract boundary control problem (\ref{eq: boundary control problem}) is well-posed~\cite[Def. 3.3.2]{Curtain2012} . In particular, we can consider the following homogeneous abstract differential equation
\begin{equation}\label{eq: boundary control problem - lifted}
\left\{\begin{split}
\dfrac{\mathrm{d} V}{\mathrm{d} t}(t) & = \mathcal{A} V(t) - B \dot{U}(t) + \mathcal{A}_d B U(t),\; t>0 \\
V(0) & = V_0 \in D(\mathcal{A})
\end{split}\right.
\end{equation}
By \cite[Th. 3.3.3]{Curtain2012}, assuming that $U \in \mathcal{C}^2(\mathbb{R}_+;\mathbb{R}^2)$, $X_0\in\mathcal{D}(\mathcal{A}_d)$ and $V_0=X_0-BU(0) \in \mathcal{D}(A)$, (\ref{eq: boundary control problem}) and (\ref{eq: boundary control problem - lifted}) each admit a unique classic solution, denoted respectively by $X(t)$ and $V(t)$, which are related by $V(t) = X(t) - B U(t)$ for all $t \geq 0$. Note that as $D(\mathcal{A}) = D(\mathcal{A}_d) \cap \mathrm{ker}(\mathcal{B})$ and $X_0,BU(0) \in D(\mathcal{A}_d)$, the condition $X_0-BU(0)\in D(\mathcal{A})$ is equivalent to $X_0-BU(0)\in \mathrm{ker}(\mathcal{B})$, i.e., based on $\mathcal{B}B=I_{\mathbb{R}^2}$, $\mathcal{B} X_0 = U(0)$. Therefore, the condition $X_0-BU(0) \in \mathcal{D}(A)$ only ensures that the boundary condition in (\ref{eq: boundary control problem}) is satisfied by the initial condition $X_0$ and the initial input $U(0)$. Let $T:\mathbb{R}_+ \rightarrow \mathcal{L}(\mathcal{H})$ be the $C_0$-semigroup generated by $\mathcal{A}$ on $\mathcal{H}$ endowed by either $\left\langle  \cdot , \cdot \right\rangle_{\mathcal{H},1}$ or $\left\langle  \cdot , \cdot \right\rangle_{\mathcal{H},2}$. Then, the unique classic solution of (\ref{eq: boundary control problem - lifted}) is given, for any $V_0 \in D(\mathcal{A})$ and all $t \geq 0$, by \cite[Th. 3.1.3]{Curtain2012}
\begin{equation*}
V(t) = T(t) V_0 + \int_0^t T(t-s) \left( - B \dot{U}(s) + \mathcal{A}_d B U(s) \right) \mathrm{d}s .
\end{equation*}
We deduce that the unique classic solution of (\ref{eq: boundary control problem}) is given, for any $X_0 \in D(\mathcal{A}_d)$ such that $\mathcal{B}X(0)=U(0)$ and all $t \geq 0$, by
\begin{align}\label{eq: boundary control problem - sol with pert}
X(t) = & T(t) (X_0 - BU(0)) + BU(t) \\
& + \int_0^t T(t-s) \left( - B \dot{U}(s) + \mathcal{A}_d B U(s) \right) \mathrm{d}s . \nonumber
\end{align}

\section{Stability Assessment}\label{sec: stability}

\subsection{Exponential stability of the $C_0$-semigroup}
In this subsection, we consider the disturbance free case, i.e. $U=0$. Then, $X(t) = T(t) X_0 \in D(\mathcal{A})$ is the unique solution of $(\mathrm{d}X/\mathrm{d}t)(t) = \mathcal{A} X(t)$ associated to the initial condition $X_0\in D(\mathcal{A})$. We define,
\begin{equation}\label{eq: definition augmented energy}
\forall t \geq 0, \; \mathcal{E}(t) \triangleq \dfrac{1}{2} \norme{X(t)}_{\mathcal{H},2}^2 = \dfrac{1}{2} \left\langle X(t) , X(t) \right\rangle_{\mathcal{H},2} .
\end{equation}
As $T(t)$ is a $C_0$-semigroup, $\mathcal{E}\in\mathcal{C}^1(\mathbb{R}_+;\mathbb{R})$ with, for any $t \geq 0$,
\begin{equation}
\dot{\mathcal{E}}(t) = \left\langle \dot{X}(t) , X(t) \right\rangle_{\mathcal{H},2} = \left\langle \mathcal{A}X(t) , X(t) \right\rangle_{\mathcal{H},2} . \label{eq: derivative energie E}
\end{equation}

Obviously, the stability properties of the $C_0$-semigroup $T(t)$ depends on the wing physical parameters. In this work, we impose the following assumptions regarding the wing parameters.
\begin{assum}\label{assumption: physical parameters}
The physical parameters involved in (\ref{eq: EDP flexion}-\ref{eq: EDP torsion}) are such that there exist $r_1,r_2,\ldots,r_8>0$ along with $0 < \epsilon_1 < \min( \epsilon_1^* , 1/K_m)$ and $0 < \epsilon_2 < \min( \epsilon_2^* , 1/K_m)$ such that
\begin{align*}
\lambda_1 \triangleq & \epsilon_1 \left( 1 - \dfrac{\sqrt{\overline{\eta_\omega}}}{2 r_1} - \dfrac{8 l^4 \overline{\rho}}{\pi^4 \underline{EI}} \left( \dfrac{\overline{|\alpha_\omega|}}{r_4} + \dfrac{\overline{|\beta_\omega|}}{r_5} + \dfrac{\overline{|\gamma_\omega|}}{\sqrt{\overline{\rho}} r_3} \right) \right) , \\  
\lambda_2 \triangleq & \overline{|\gamma_\omega|} + \dfrac{\sqrt{\overline{\rho}\underline{\rho}} \, \overline{|\alpha_\omega|} + \epsilon_2 \overline{I_w}\,\overline{|\gamma_\phi|}}{2 \sqrt{\underline{\rho}} r_6} + \dfrac{\epsilon_1 \sqrt{\overline{\rho}} \, \overline{|\gamma_\omega|} r_3}{2} \\
& + \dfrac{\sqrt{\overline{\rho}}\,\overline{|\beta_\omega|}/\sqrt{\underline{I_w}} + \sqrt{\overline{I_w}}\,\overline{|\gamma_\phi|}/\sqrt{\underline{\rho}}}{2 r_7} , \\  
\lambda_3 \triangleq & 1 - \epsilon_1 \left( \dfrac{16 l^4 \overline{\rho}}{\pi^4 \underline{\eta_\omega EI}} + \dfrac{\sqrt{\overline{\eta_\omega}} r_1}{2} \right) , \\  
\lambda_4 \triangleq & \epsilon_2 \left( 1 - \dfrac{\sqrt{\overline{\eta_\phi}}}{2 r_2} \right) - \dfrac{4 l^2}{\pi^2 \underline{GJ}} \left( \dfrac{(\sqrt{\overline{\rho}\underline{\rho}}\,\overline{|\alpha_\omega|} + \epsilon_2 \overline{I_w} \overline{|\gamma_\phi|})r_6}{2 \sqrt{\underline{\rho}}} \right. \\
& \left. + \dfrac{\sqrt{\overline{I_w}}(\overline{|\alpha_\phi|} + \epsilon_2 \overline{|\beta_\phi|})}{2 r_8} + \dfrac{\epsilon_1 \overline{\rho}\,\overline{|\alpha_\omega|} r_4}{2} + \epsilon_2 \overline{I_w}\,\overline{|\alpha_\phi|} \right) , \\  
\lambda_5 \triangleq & \overline{|\beta_\phi|} + \dfrac{(\sqrt{\overline{\rho}}\,\overline{|\beta_\omega|}/\sqrt{\underline{I_w}} + \sqrt{\overline{I_w}}\,\overline{|\gamma_\phi|}/\sqrt{\underline{\rho}})r_7}{2} \\
& + \dfrac{\sqrt{\overline{I_w}} (\overline{|\alpha_\phi|} + \epsilon_2 \overline{|\beta_\phi|})r_8}{2} + \dfrac{\epsilon_1 \overline{\rho} \overline{|\beta_\omega|} r_5}{2 \underline{I_w}} , \\  
\lambda_6 \triangleq & 1 - \epsilon_2 \left( \dfrac{4 l^2 \overline{I_w}}{\pi^2 \underline{\eta_\phi GJ}} + \dfrac{\sqrt{\overline{\eta_\phi}} r_2}{2} \right) , \\  
\end{align*}
satisfy $\lambda_1, \ldots, \lambda_6 > 0$, $\pi^4 \underline{\eta_\omega EI} \lambda_3 / (16 l^4 \overline{\rho}) - \lambda_2 > 0$, and $\pi^2 \underline{\eta_\phi GJ} \lambda_6 / (4l^2 \overline{I_w}) - \lambda_5 > 0$.
\end{assum}

As with the constraints imposed in~\cite{Bialy2016,He2017}, Assumption~\ref{assumption: physical parameters} imposes a trade-off between the structural stiffness of the wing and the amplitude of the aerodynamic coefficients. Indeed, it is easy to see that for fixed aerodynamic coefficients, Assumption~\ref{assumption: physical parameters} asymptotically boils down to the positiveness constraints of the coefficients of (\ref{eq: inner A1X-X}), for which it has been shown that a feasible solution always exists, when the stiffness increases. Therefore, Assumption~\ref{assumption: physical parameters} can always be satisfied by adequately increasing the stiffness of the structure.

Under Assumption~\ref{assumption: physical parameters}, we can now assess the exponential stability of the closed-loop system.
\begin{theorem}\label{th: exponential stab C0-semigroup 2nd inner product}
Provided that Assumption~\ref{assumption: physical parameters} holds, $T(t)$ is an exponentially stable $C_0$-semigroup.
\end{theorem}
\textbf{Proof.}
Let $X(t) = T(t)X_0 = (f(\cdot,t),g(\cdot,t),h(\cdot,t),z(\cdot,t))\in D(\mathcal{A})$. Based on the upper-bound of $\left\langle \mathcal{A}_1 X(t) , X(t) \right\rangle_{\mathcal{H},2}$ given in (\ref{eq: inner A1X-X}) and applying Young's and Poincar{\'e}'s inequalities to the upper bound of $\left\langle \mathcal{A}_2 X(t) , X(t) \right\rangle_{\mathcal{H},2}$, it yields for any $t\geq0$,
\begin{align}
\dot{\mathcal{E}}(t) 
\leq & - k_1 (g(l,t)+\epsilon_1 f(l,t))^2 - k_2 (z(l,t)+\epsilon_2 h(l,t))^2 \nonumber \\
& - \lambda_1 \int_0^l EI(y) f''(y,t)^2 \mathrm{d}y - \lambda_4 \int_0^l GJ(y) h'(y,t)^2 \mathrm{d}y \nonumber \\
& - \left( \dfrac{\pi^4 \underline{\eta_\omega EI}\lambda_3}{16 l^4 \overline{\rho}} - \lambda_2 \right)  \int_0^{l} \rho(y) g(y,t)^2 \mathrm{d}y \nonumber \\
& - \left( \dfrac{\pi^2 \underline{\eta_\phi GJ} \lambda_6}{4 l^2 \overline{I_w}} - \lambda_5 \right) \int_0^{l} I_w(y) z(y,t)^2 \mathrm{d}y . \label{eq: time derivative of the augmented energy}
\end{align}
with $\lambda_1,\ldots,\lambda_6>0$ defined in Assumption~\ref{assumption: physical parameters}. Introducing
\begin{equation*}
\mu_m \triangleq 2 \min\left( \lambda_1 , \dfrac{\pi^4 \underline{\eta_\omega EI}\lambda_3}{16 l^4 \overline{\rho}} - \lambda_2 , \lambda_4, \dfrac{\pi^2 \underline{\eta_\phi GJ} \lambda_6}{4 l^2 \overline{I_w}} - \lambda_5 \right) > 0 ,
\end{equation*}
and recalling that $k_1,k_2 \geq 0$, we obtain, based on the equivalence of the norms (\ref{eq: norm equivalent}),
\begin{equation*}
\forall t\geq0 , \; \dot{\mathcal{E}}(t) \leq - \dfrac{\mu_m}{2} \norme{X(t)}_{\mathcal{H},1}^2 \leq  - \Lambda \mathcal{E}(t) ,
\end{equation*}
where $\Lambda \triangleq \mu_m/(1+\epsilon_m K_m)>0$ is independent of the initial condition $X_0\in D(\mathcal{A})$. Thus, as $\mathcal{E}\in\mathcal{C}^1(\mathbb{R}_+;\mathbb{R})$, we get that for any $t \geq 0$, $\mathcal{E}(t) \leq \mathcal{E}(0) e^{-\Lambda t}$. Furthermore, as by definition $\mathcal{E}(t)=\norme{X(t)}_{\mathcal{H},2}^2/2$, it yields that
\begin{equation*}
\forall X_0 \in D(\mathcal{A}),\; \forall t\geq0 ,\; \norme{T(t)X_0}_{\mathcal{H},2} \leq \norme{X_0}_{\mathcal{H},2} e^{-\Lambda t/2} . 
\end{equation*}
As $T(t)\in\mathcal{L}(\mathcal{H})$ for all $t\geq0$ and $\overline{D(\mathcal{A})} = \mathcal{H}$, the above inequality can be extended for all $X_0 \in \mathcal{H}$. Thus, $T(t)$ is an exponentially stable $C_0$-semigroup for $\norme{\cdot}_{\mathcal{H},2}$ with $\left\Vert T(t)\right\Vert_{\mathcal{H},2} \leq e^{-\Lambda t/2}$ for all $t\geq0$. Resorting to (\ref{eq: norm equivalent}), it shows that $\left\Vert T(t)\right\Vert_{\mathcal{H},1} \leq K_E e^{-\Lambda t/2}$ for all $t\geq0$ where $K_E = \sqrt{(1+K_m \epsilon_m)/(1-K_m \epsilon_m)}$. In particular, the growth bound $\omega_0(T)$ of $T(t)$ is such that $\omega_0(T) \leq -\Lambda/2 < 0$.
\qed

\subsection{Stability analysis for bounded input perturbations}

Let $X(t)$ be a solution of the boundary control problem (\ref{eq: boundary control problem}). Under Assumption~\ref{assumption: physical parameters}, it can be obtained by (\ref{eq: boundary control problem - sol with pert}) that for all $t\geq 0$,
\begin{align}
\left\Vert X(t) \right\Vert_{\mathcal{H},1} \leq & K_E \left\Vert X_0 - BU(0) \right\Vert_{\mathcal{H},1} e^{- \Lambda t /2} + \left\Vert B \right\Vert \left\Vert U(t) \right\Vert_2 \nonumber\\
& + K_E \left\Vert \mathcal{A}_d B \right\Vert \int_0^t e^{-\Lambda(t-s)/2} \left\Vert U(s) \right\Vert_2 \mathrm{d}s \nonumber \\
& + K_E \left\Vert B \right\Vert \int_0^t e^{-\Lambda(t-s)/2} \left\Vert \dot{U}(s) \right\Vert_2 \mathrm{d}s . \label{eq: upper bound norm state with peturbations}
\end{align}

\subsubsection{Bounded input disturbances}
Assume that $U$ and $\dot{U}$ are bounded. It yields,
\begin{align*}
\left\Vert X(t) \right\Vert_{\mathcal{H},1} \leq & K_E \left\Vert X_0 - BU(0) \right\Vert_{\mathcal{H},1} e^{- \Lambda t /2} \\
& + \left( \left\Vert B \right\Vert + \dfrac{2 K_E}{\Lambda} \left\Vert \mathcal{A}_d B \right\Vert \right) \underset{s\in\mathbb{R}_+}{\sup} \left\Vert U(s) \right\Vert_2 \\
& + \dfrac{2 K_E}{\Lambda} \left\Vert B \right\Vert \underset{s\in\mathbb{R}_+}{\sup} \left\Vert \dot{U}(s) \right\Vert_2 .
\end{align*}
Thus, the system energy is bounded and, as $t$ tends to infinity, we have
\begin{align*}
\underset{t \rightarrow +\infty}{\mathrm{lim\,sup}} \left\Vert X(t) \right\Vert_{\mathcal{H},1} \leq & \left( \left\Vert B \right\Vert + \dfrac{2 K_E}{\Lambda} \left\Vert \mathcal{A}_d B \right\Vert \right) \underset{s\in\mathbb{R}_+}{\sup} \left\Vert U(s) \right\Vert_2 \\
& + \dfrac{2 K_E}{\Lambda} \left\Vert B \right\Vert \underset{s\in\mathbb{R}_+}{\sup} \left\Vert \dot{U}(s) \right\Vert_2 .
\end{align*}
In particular, the contribution of the initial condition vanishes exponentially. Employing Agmon's and then Poincar{\'e}'s inequalities yields 
\begin{equation*}
\left\Vert h(\cdot,t) \right\Vert_\infty^4 
\leq \dfrac{16 l^2}{\pi^2} \left\Vert h'(\cdot,t) \right\Vert_{L^2(0,l)}^4
\leq \dfrac{16 l^2}{\pi^2 \underline{GJ}^2} \left\Vert X(t) \right\Vert_{\mathcal{H},1}^4 .
\end{equation*}
Applying a similar procedure to $f$ and $f'$, it shows that in the disturbance free case (i.e., $U=0$), both bending and twisting displacements converge exponentially and uniformly over the wingspan to zero. In the presence of bounded input disturbances, the contribution of the initial condition to the displacements vanishes exponentially. Furthermore, the displacements are bounded in time, uniformly over the wingspan, and 
\begin{align*}
\underset{t \rightarrow +\infty}{\mathrm{lim\;sup}} \left|\left|f(\cdot,t)\right|\right|_\infty & \leq \dfrac{4 l^{3/2}}{\pi^{3/2} \underline{EI}^{1/2}} \; \underset{t \rightarrow +\infty}{\mathrm{lim\;sup}} \left\Vert X(t) \right\Vert_{\mathcal{H},1}
< +\infty , \\
\underset{t \rightarrow +\infty}{\mathrm{lim\;sup}} \left|\left|f'(\cdot,t)\right|\right|_\infty & \leq \dfrac{2 l^{1/2}}{\pi^{1/2} \underline{EI}^{1/2}} \; \underset{t \rightarrow +\infty}{\mathrm{lim\;sup}} \left\Vert X(t) \right\Vert_{\mathcal{H},1}
< +\infty , \\
\underset{t \rightarrow +\infty}{\mathrm{lim\;sup}} \left|\left|h(\cdot,t)\right|\right|_\infty & \leq \dfrac{2 l^{1/2}}{\pi^{1/2} \underline{GJ}^{1/2}} \; \underset{t \rightarrow +\infty}{\mathrm{lim\;sup}} \left\Vert X(t) \right\Vert_{\mathcal{H},1}
< +\infty .
\end{align*}
The obtained upper-bounds on the system energy and on both bending and twisting displacements are function of $\left\Vert B \right\Vert$ and $\left\Vert \mathcal{A}_d B \right\Vert$ given by (\ref{eq: norm B}) and (\ref{eq: norm UB}), respectively. It ensures that the increase of the wing stiffness will reduce the impact of the perturbations on the closed-loop system.

\subsubsection{Vanishing input disturbances}
Assume that the disturbance input is vanishing in the following sense:
\begin{equation}\label{eq: vanishing input}
\underset{t \rightarrow + \infty}{\lim} \left\Vert U(t) \right\Vert_2 = \underset{t \rightarrow + \infty}{\lim} \left\Vert \dot{U}(t) \right\Vert_2 = 0 .
\end{equation}
Then, $\left\Vert X(t) \right\Vert_{\mathcal{H},1}$ converges to zero as $t$ tends to infinity. Indeed, based on (\ref{eq: upper bound norm state with peturbations}), this claim will be true if the two integral terms converge to zero as $t$ tends to infinity. Consider an arbitrary $\epsilon > 0$. By (\ref{eq: vanishing input}), there exists $T \geq 0$ such that for any $t \geq T$, $\left\Vert U(t) \right\Vert_2 \leq \Lambda \epsilon /2$. Fixing such a $T \geq 0$, we have for all $t \geq T$,
\begin{align*}
& \int_0^t e^{-\Lambda(t-s)/2} \left\Vert U(s) \right\Vert_2 \mathrm{d}s \\ 
= & e^{-\Lambda t/2} \int_0^T e^{\Lambda s/2} \left\Vert U(s) \right\Vert_2 \mathrm{d}s 
+ \int_T^t e^{-\Lambda(t-s)/2} \left\Vert U(s) \right\Vert_2 \mathrm{d}s \\
\leq & e^{-\Lambda t/2} \int_0^T e^{\Lambda s/2} \left\Vert U(s) \right\Vert_2 \mathrm{d}s 
+  \epsilon ,
\end{align*}
which yields
\begin{equation*}
\underset{t \rightarrow + \infty}{\mathrm{lim\;sup}} \int_0^t e^{-\Lambda(t-s)/2} \left\Vert U(s) \right\Vert_2 \mathrm{d}s \leq \epsilon .
\end{equation*}
As this inequality is true for any $\epsilon > 0$, and due to the fact that the integral is positive for all $t \geq 0$, it implies that
\begin{equation*}
\underset{t \rightarrow + \infty}{\lim} \int_0^t e^{-\Lambda(t-s)/2} \left\Vert U(s) \right\Vert_2 \mathrm{d}s = 0 .
\end{equation*}
Similarly, the second integral of (\ref{eq: upper bound norm state with peturbations}) converges to zero when $t \rightarrow \infty$. Thus we have $\left\Vert X(t) \right\Vert_{\mathcal{H},1} \underset{t \rightarrow +\infty}{\longrightarrow} 0$. Then, we deduce from the bounded case that 
\begin{equation*}
\underset{t \rightarrow + \infty}{\lim} \left|\left|f(\cdot,t)\right|\right|_\infty 
= \underset{t \rightarrow + \infty}{\lim} \left|\left|f'(\cdot,t)\right|\right|_\infty 
= \underset{t \rightarrow + \infty}{\lim} \left|\left|h(\cdot,t)\right|\right|_\infty
= 0 ,
\end{equation*}
i.e., both bending and twisting displacements converge uniformly over the wingspan to zero when $t$ tends to infinity.

\section{Numerical Simulations}\label{sec: simulations}
\begin{figure}[htb]
	\centering
	\includegraphics[width=3in]{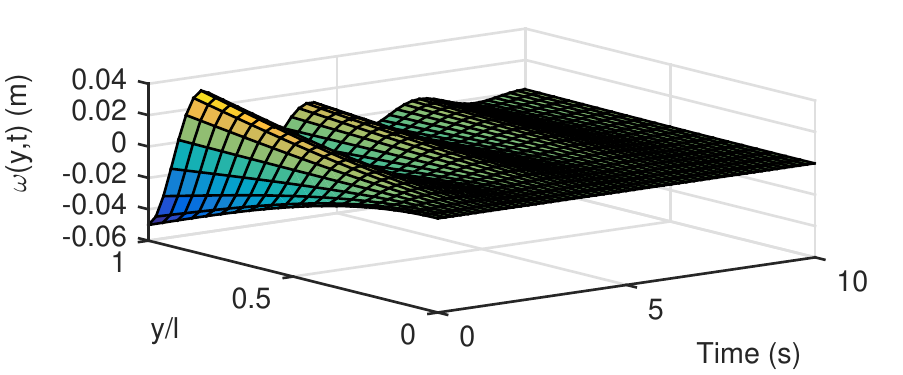}
	\includegraphics[width=3in]{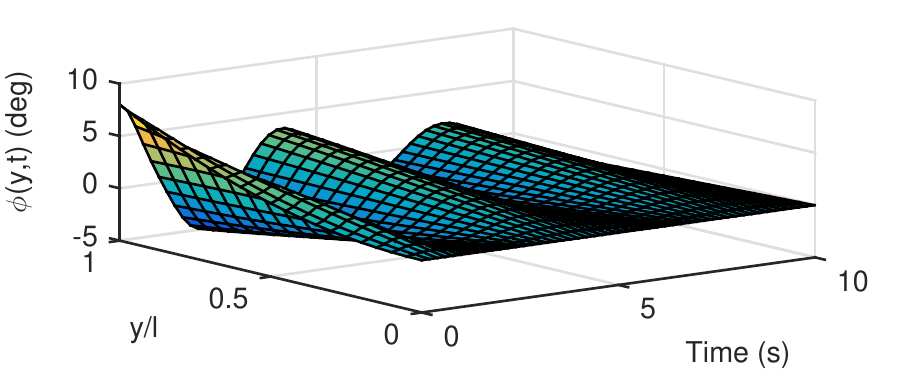}
	\caption{Open-loop response}
	\label{fig: OL}
\end{figure}

\begin{figure}[htb]
	\centering
	\includegraphics[width=3in]{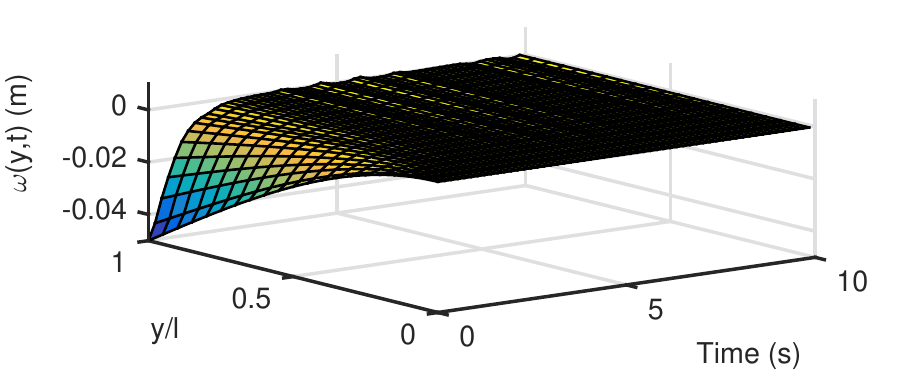}
	\includegraphics[width=3in]{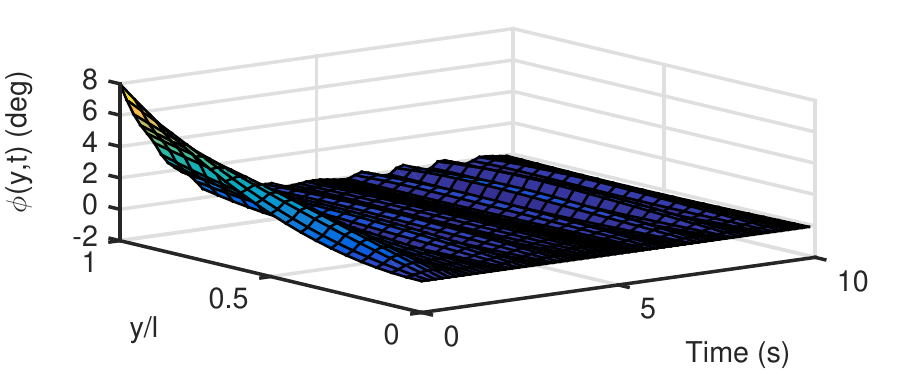}
	\caption{Closed-loop response in the presence of input perturbations}
	\label{fig: CL perturbations}
\end{figure}

The numerical scheme is based on the Galerkin method \cite{ciarlet2002finite}. For simulations purposes, the following persistent input perturbations are considered.
\begin{align*}
u_1(t) & = 3 \cos(0.2 \pi t) \sin(\pi t) \cos(3 \pi t) , \\
u_2(t) & = \sin(0.2 \pi t) \cos(\pi t) \sin(3 \pi t) .
\end{align*}
The initial condition is selected as $\omega_0(y) = y^2(y-3l)/(40 l^2)$, $\omega_{t0}(y)=0$, $\phi_0(y) = 2 \pi y^2 / (45 l^2)$ and $\phi_{t0}(y) = 0$. The open-loop response is depicted in Fig.~\ref{fig: OL}, which exhibits poorly damped oscillations. Setting the controller gains as $k_1 = 10$ and $k_2 = 4$, the behavior of the closed-loop system is shown in Fig.~\ref{fig: CL perturbations}. It can be seen that the flexible displacements are damped out rapidly, even in the presence of input perturbations.

\section{Conclusion}\label{sec: conclusion}
The well-posedness and the stability properties of a flexible nonhomogeneous wing in the presence of input disturbances have been studied. The wing is modeled by a distributed parameter system for which the well-posedness issue has been tackled in the framework of semigroups. The stability of the closed-loop system has been investigated by a Lyapunov-based approach. It has been shown that, under physical structural constraints, both flexible displacements are bounded and will exponentialy converge to zero for vanishing disturbances.


%



\ifCLASSOPTIONcaptionsoff
  \newpage
\fi



\bibliographystyle{IEEEtranS}
\nocite{*}
\bibliography{IEEEabrv,mybibfile}

\end{document}